# Logical Primes, Metavariables and Satisfiability

Bernd R. Schuh

Dr. Bernd Schuh, Bernhardstraße 165, D-50968 Köln, Germany; bernd.schuh@netcologne.de

*Abstract*. For formulas F of propositional calculus I introduce a "metavariable" $M_F$ and show how it can be used to define an algorithm for testing satisfiability. $M_F$ is a formula which is true/false under all possible truth assignments iff F is satisfiable/unsatisfiable. In this sense $M_F$ is a metavariable with the "meaning" 'F is SAT'. For constructing $M_F$ a group of transformations of the basic variables $a_i$ is used which corresponds to 'flipping" literals to their negation. The whole procedure corresponds to branching algorithms where a formula is split with respect to the truth values of its variables, one by one. Each branching step corresponds to an approximation to the metatheorem which doubles the chance to find a satisfying truth assignment but also doubles the length of the formulas to be tested, in principle. Simplifications arise by additional length reductions.
I also discuss the notion of "logical primes" and show that each formula can be written as a uniquely defined product of such prime factors. Satisfying truth assignments can be found by determining the "missing" primes in the factorization of a formula.

*Introduction*

Introductions to the problem of satisfiability can be found in textbooks and reviews, some of them available in the net (see e.g. [1],[2]). One of the unsolved questions of the field is whether satisfiability can be determined in polynomial time ("P=NP ?"). Other questions center around efficient techniques to determine satisfying assignments (see [3,4] for new approaches), and to identify classes of "hard" problems which inherently seem to consume large computing time. I believe that some insight into the difficulties can be gained by using algebraic tools. I have outlined some of them in a previous note [5]. In particular the notion of 'logical primes' and the group of flipping transformations appear helpful in analyzing formulas and deriving general theorems.

I will recall these notions and some consequences in sections I , II and III. Then I will introduce the metaformula and a related quantity, the parityformula, which encodes whether F has an even or an odd number of satisfying solutions in section IV. In sections V and VI algorithms with encreasing effectiveness in determining satisfiability are introduced.

*I. Definitions*

We consider a finite algebra V with two operations + and x, and denote by $\mathbf{1}$ and $\mathbf{0}$ their neutral elements, respectively, i.e.

(1)   $a \times \mathbf{1} = a$,  $a + \mathbf{0} = a$

Additionally, the operations are associative and commutative, and the distributive law

(2)   $a \times (b+c) = a \times b + a \times c$

is assumed to hold in V.

Two more properties are required, namely:

(3)   $a+a=0$

(4)   $a \times a = a$

It is clear from these definitions that V may be identified with the Boolean algebra of propositional calculus, where "x" corresponds to the logical "AND" and "+" to the logical "XOR" (exclusice OR).

To each element of V we introduce its "negation" by

(5)   $\sim a := a+1$

From (2), (3) and (4) it is clear that $\sim a \times a = 0$ as is appropriate for a negation.

II. *Consequences.*

As a first consequence of equ.s (1) - (5) we can state the following theorem:

(TI)   $\dim(V) = |V| = 2^N$   for some natural number N

i.e. the number of elements of V is necessarily a power of 2.

This is not surprising, of course, if one has the close resemblence of V to propositional calculus in mind. But here it is to be deduced solely from the algebraic properties.

All proofs are given in the appendix.

In order to formulate a second consequence it is necessary to introduce the notion of "logical primes". We define:

(DI)   $p \in V$ is a (logical) prime, iff

for any $a \in V$   $p \times a = 0$ implies $a = 0$ or $a = \sim p$.

If not clear by definition, the name "prime" will become clear by the following theorems

(TII)   There are exactly $ld|V|=N$ many primes in V. And:

(TIII)  Each element of V has a unique decomposition into primes:

(6)   $a = \Pi_j p_j$   where the product refers to the x-operation, and $j \in I_a$, and $I_a = I_b$ iff $a=b$

This property can be formulated alternatively with the negated primes $\sim p_j$ via

(7)   $a = \Sigma_j \sim p_j$ with $j \in {}^c I_a$  (${}^c I_a$ is the complement of $I_a$ in $\{0,1,..., N-1\}$ )

The neutral elements $0$ and $1$ are special cases. $1$ is expressed as the empty product according to (6), whereas the sum extends over all primes. For $0$ the sum-representation is empty, but the product extends over all possible primes.

A property which is extremely helpful in calculations is

(8)   $\sim p_j \times \sim p_k = \sim p_k \, \delta_{jk}$   ($\delta_{jk} = 1$ iff $j=k$, 0 otherwise)

which with the aid of (5) can be written

$p_j \times p_k = p_j + \sim p_k = \sim p_j + p_k$   for $k \neq j$

Note, that no use has been made of the correspondence of $\{V,+,\times, 0, 1\}$ to propositional calculus, up to now. We can even proceed further and define the analogue of truth assignments. Consider the set of maps $T: V \rightarrow \{0,1\}$. We call T

"allowed" iff there is a relationship between the image of a "sum" or a "product" and the image of the single summands or factors. In formula:

(9)     $T(a+b) = f(T(a),T(b))$  and  $T(a \times b) = g(T(a),T(b))$

with some functions f and g and all $a,b \epsilon V$.

These relations suffice to show theorem IV

(TIV)  There are exactly N different allowed maps $T_j$, and they fulfill:

(10)    $T_j(\sim p_k) = \delta_{jk}$

Given functions f and g of (9) one can also use (10) as a definition and extend $T_j$ to all elements of V via (7).

In one last step we assume $N=2^n$ for some natural number n. Then

(TV)   n distinct elements $a_k$ ( different from $0, 1$) can be found, such that

(11)    $\sim p_s = (\Pi_j s_j a_j)(\Pi_k (1-s_k) \sim a_k)$    where $s = \Sigma_r 2^{r-1} s_r$ is the binary representation of s.

In words: each element of V can be written as a "sum" of "products" of all $a_k$ and $\sim a_k$. E.g. for n=3 one has $p_2 = a_2 \times \sim a_1 \times \sim a_3$ as one of the eight primes. The $a_k$ are not necessarily unique. E.g., for n=3, given $a_k$, the set $a_1, a_3, a_1 \times \sim a_2 + \sim a_1 \times a_2$ will serve the same purpose (with a different numbering convention in (11)).

*III. Propositional calculus.*

Propositional calculus (PC) consists of infinitely many formulas which can be constructed from basic variables $a_k$ with logical functions (like "AND", "OR" and negation). Even for a finite set of n basic variables $B_n=\{a_1,a_2,...a_n\}$ there are infinitely many formulas arizing from combinations of the basic variables. These formulas can be grouped into classes of logically equivalent formulas. That is, formulas F and F' belong to the same class iff their values under any truth assignment $\mathcal{T}:B_n \to \{0,1\}$ are the same. Members of different classes are logically inequivalent, i.e. there is at least one truth assignment for which their values differ. This finite set of classes for fixed n can be identified with the algebra V of the foregoing section. Neutral elements of the operations x and +, $\mathbf{1}$ and $\mathbf{0}$, are interpreted as complete truth and complete unsatisfiability.

In order to see how operations + and x correspond to logical operations "AND" and "OR" we define a new operation v in V via

(12)    a v b = a + b + axb

With this definition the defining relations (1) - (5) can be reformulated in terms of v and x, and the algebraic structure of a Boolean algebra for formulas becomes obvious. v is the logical "OR", x the logical "AND".

Relation (12) reduces logical considerations to simple algebraic manipulations in which + and x can be used as in multiplication and addition of numbers. Additionally the simplifying relations a+a=$\mathbf{0}$ = ~axa  and axa=a, a+~a=$\mathbf{1}$  hold. Consider for illustration  the so called "resolution" method. It states that avb and ~avc imply bvc. A "calculational" proof of this statement might run as follows.  (From now on we skip the x-symbol for multiplication) We make use of the fact that in PC  the implication a$\Rightarrow$b is identical to  ~a v b :

(avb)(~avc) ⟹ bvc = ~((avb)(~avc))vbvc =

($\mathit{1}$+(a+~ab)(~a+ca))+b+c+bc+(b+c+bc)($\mathit{1}$+(a+~ab)(~a+ca)) = $\mathit{1}$ +ac+~ab +

+(b+c+bc)(a+~ab)(~a+ca) = $\mathit{1}$ +ac+~ab+abc+~ab+~bac = $\mathit{1}$ +ac($\mathit{1}$ +b+~b) = $\mathit{1}$

In other words: the implication is a tautology ( true under all truth assignments) as claimed.

TIII and TV tell us that each formula F of PC has a unique decomposition into a "sum" of "products" of its independent variables $a_k$. Because of (8) and (12) the sum in (7) may be written as a "v"-sum. Thus (7) takes the form of a disjunctive normal form (DNF) and it can as well be transformed into a conjunctive normal form (CNF) as given by (6). For the neutral element $\mathit{0}$ one has

(13)    $\mathit{0}$ = ($a_1$v$a_2$v...v$a_n$)x(~$a_1$v...$a_n$)x...x(~$a_1$v...~$a_n$)

with all possible primes. According to (6) each formula F has a similar representation, but with some prime factors missing. From the primes present one can immediately read off the truth assignments for which F evaluates to 0, thus the missing factors give the truth assignments for which F is satisfiable.

Note, however, that each factor in the prime representation of a formula involves *all* $a_k$ . So one way of determining satisfying assignments or test a formula for satisfiability consists of transforming a given CNF representation of the formula to its standard form (6). This can be done e.g. by "blowing up" each factor until all $a_k$ are present. E.g.   avbv~c = (avbv~cvd)(avbv~cv~d)  from 3 to 4 variables. Since each new factor has to be treated in the same way, until n is reached, this is a O($2^n$) - process in principle, which makes the difficulty in finding a polynomial time algorithm for testing satisfiability understandable.

Also from (7) with (10) and (8) it follows that the satisfying assignments of a formula $F= \Sigma_j \sim p_j$ are given by the negated primes which do not show up in the CNF representation. In particular, the number of satisfying assignments is equal to the number of summands in this equation. Furthermore, they can be read off immediately, since, according to (10) $T_s(F) = 1$ iff the corresponding $\sim p_s$ shows up in the sum. Also the $T_s$ must coincide with the $2^n$ possible truth assignments $\mathcal{T}:B_n \rightarrow \{0,1\}$. One may choose the numbering such that the values of $T_s$ on $B_n$ are given by the binary representation $s= \Sigma_r 2^{r-1} T_s(a_r)$.

As a last example for the usefulness of the algebraic approach we consider the number of satisfying assignments of a formula F of PC , #(F) and show that this number does not change if some (or all) of the variables $a_k$ are "flipped", i.e. substituted by their negation and vice versa:

(14)   $\#(F(a_1,...,a_n)) = \#(F(a_1,...\sim a_i,...\sim a_j,...))$

To prove this "conservation of satisfiability" we consider a group of transformations $\{R_0,...R_{N-1}\}$ which negate the $a_k$ according to the following definition: $R_s$ negates all $a_r$ (and $\sim a_r$ likewise) for which $s_r$ in the binary representation of s is non zero. In formula, for any truth assignment $T_j$

(15)   $T_j(R_s(a_r)) = (1-s_r) T_j(a_r) + s_r(1- T_j(a_r))$    and    $s= \Sigma_r 2^{r-1} s_r$.

It is easy to see that the $R_s$ form a group with $R_0 = \text{id}$, and each $R_s$ induces a permutation $\pi_s$ of of the $\tilde{p}_j$ which is actually a transposition given by

(15')  $\pi_s(j) = s + j - 2\Sigma_r 2^{r-1} s_r j_r =: (s,j)$

Thus $R_s$ simply permutes the primes $p_k$ and therefore in the representation of F in (6) or (7) their number is not changed. The fact may also be stated as

(16)  $T_j(R_s(F)) = T_{(s,j)}(F)$

and therefore $\#(F) = \Sigma_j T_j(F) = \Sigma_j T_j(R_s(F)) = \#(R_s(F))$ which proves (14).

One may also conclude from (16) that for satisfiable F each $R_s(F)$ is satisfiable. More precise: if $T_j(F) = 1$ for some j, then for any k there is a flipping operation $R_s$ such that $T_k(R_s(F)) = 1$, namely $s = (k,j)$. Likewise, for any $R_s$ one can find a $T_k$ such that $T_k(R_s(F)) = 1$.

On the other hand, if F is not satisfiable, none of the $R_s(F)$ can be satisfiable, otherwise one would have $T_k(R_s(F)) = 1$ for some k and thus $T_j(F) = 1$ for some j, contrary to the assumption that F is not SAT.

*IV. The Metaformula.*

For any formula F of n variables we write $F(a_1,\ldots,a_n)$ and define the metaformula by "adding" with respect to the OR-operation all "flipped" versions of F:

(17)  $M_F(a_1,\ldots,a_n) = R_0(F) \vee R_1(F) \vee \ldots \vee R_{N-1}(F)$

where $N = 2^n$. From the considerations at the end of the foregoing section it is immediately clear that $M_F$ is not satisfiable if F is not, and that $M_F$ is a tautology if F is satisfiable:

(18)   $M_F = 1$  iff  $F \varepsilon$ SAT;   $M_F = 0$  iff  $F \not\varepsilon$ SAT

Thus, considered as a logical variable itself, $M_F$ represents the satisfiability of F. $M_F$ can only take the two values 0 and 1 depending on whether F is SAT or not. That is why I call $M_F$ a metatheorem or metaformula.

Very similarly one may introduce a "parityformula" $P_F$ in substituting the OR-operation in the definition (17) by the exclusive XOR. Analogously to (18) one can show that

$P_F = 0$   iff  $P_F$ has an even number of satisfying truth assignments,

$P_F = 1$   iff $P_F$ has un odd number of satisfying assignments.

*V. SAT algorithm.*

We now turn to the question how $M_F$ can be utilized to formulate SAT algorithms. Since either  $M_F = 1$ or  $M_F = 0$ it is sufficient to test one single truth assignment in order to determine whether F is SAT or not. Thus the satisfiability of F can be determined in linear time in the length of $M_F$. Nothing is gained so far, however, since the length of $M_F$ is of order N times the length of F. Thus, instead of testing all N $T_j$ on F to determine its satisfiability in the metatheorem approach one first constructs an order-N variant of F and checks it with a single $T_j$.

Simplifications may arise, however in the process of constructing $M_F$.

Note first that $M_F$ can be constructed in n steps. For this purpose consider the shift operator

(19) $\quad D^{(k)}(F) = F \vee R_q(F)$ with $q=2^{k-1}$ and $k=1,...,n$

Note that all N operators $R_s$ can be generated by the n operators $R_q$ with $q=2^{k-1}$ and $k=1,...,n$. E.g. $R_{29} = R_{16}R_8R_4R_1$.

Furthermore it is easy to see that $R_q$ flips the variable $a_k$ and therefore $D^{(k)}$ is independent of $a_k$, and (19) may be rewritten as

(20) $\quad D^{(k)}(F) = F(a_1,..., a_k= \mathbf{1}, ...,a_n) \vee F(a_1,..., a_k= \mathbf{0}, ..., a_n)$.

In terms of shift operators $M_F$ may be rewritten as

(21) $\quad M_F = D^n(D^{n-1}(...D^2(D^1(F)...) = :D^{(n)}...D^{(1)}(F)$

We can now consider systematic approximations on $M_F$. Namely the series of $l^{th}$ order approximations

(22) $\quad M_F^{(l)} = D^{(l)}...D^{(1)}(F) = D^{(l)}(M_F^{(l-1)})$ ; $M_F^{(0)} = F$.

From this definition we may write $M_F^{(l)}$ as

(23) $\quad M_F^{(l)} = F \vee R_1(F) \vee ... \vee R_{q-1}(F)$ with $q=2^l$.

We will show next that a properly chosen truth assignment for testing the l-th approximation can give a wealth of information. Check $T_q$ with $q=2^l$ on $M_F^{(l)}$. Let us assume that $T_q(M_F^{(l)}) = 1$. Then one of the $R_i(F)$ is true under that truth assignment, therefore there is a truth assignment which satisfies F, thus F is satisfiable. If on the other hand $T_q(M_F^{(l)}) = 0$, then we may conclude the following: $T_q(R_i(F)) = 0$ for all i $\varepsilon$ {1, 2,..., q}. Therefore also $T_{(i,q)}(F) = 0$ according to (16). In conclusion, see (15') for the definition of (i,q):

(24)  If $T_q(M_F^{(l)}) = 0$ then F is <u>not</u> satisfied by truth assignments $T_k$ for

  k $\varepsilon$ {q, q+1,..., 2q-1}   (q= $2^l$)

An effective check for satisfiability of F may therefore run as follows:

    CHECKSAT [F,n]

    set s=1

1    set F = $D^{(s)}(F)$

    if $T_s(F)=1$  then stop and return "F is SAT"

    s=s+1

    if s=n  then  stop and return "F is not SAT"

    goto 1

CHECKSAT determines satisfiability in n steps each of which is linear in the length of the formula. In each step the number of excluded truth assignments is doubled, as well as the chance to find a satisfying assignment if there is one. However, the formulas to be checked become longer and longer in each step, therefore it remains an order-N process in principle. A look at (20) reveals that the procedure corresponds to a successive elimination of variables. Further optimizations require length reductions in the formulas $M_F^{(l)}$ which arise in the approximation process.

*VI. Length reduction.*

In this section we assume F to be given in conjunctive normal form (CNF):

(25)     $F = C_1 C_2 ... C_m$

with m clauses of the form

(26)     $C = L \vee R$

where L is a literal corresponding to one of the variables $a_k$ or its negation, and R may itself be written in the form (26) and so forth until R is a literal.

In the process of eliminating variables described in the foregoing section the following well known rules can help to reduce formulas in length.

(27)
- (a)    $(L \vee R)(\sim L \vee R) = R$
- (b)    $(L \vee R)(\sim L \vee R') = LR' \vee \sim LR$
- (c)    $L(\sim L \vee R) = LR$
- (d)    $(L \vee R)(L \vee R') = L \vee \sim LRR'$
- (e)    $L(L \vee R) = L$

In a CNF-formula one encounters terms of the form $(L \vee R_1)...(L \vee R_s)(\sim L \vee S_1)...(\sim L \vee S_t)$ which may be rewritten by the aid of (27):

(28)     $F(L) := (L \vee R_1)...(L \vee R_s)(\sim L \vee S_1)...(\sim L \vee S_t) = L S_1 S_2 ... S_t \vee \sim L R_1 R_2 ... R_s$

Note that the CNF form on the l.h.s. is split into a disjunction of two CNF-formulas on the r.h.s.. The variable l does not show up in the R an S by definition. If we eliminate the variable L, which is exactly what happens when the shift operator D is applied, one reads off the r.h.s.:

(29)     $F(L) \vee F(\sim L) = S_1 S_2 ... S_t \vee R_1 R_2 ... R_s$

a disjunction of CNF-forms independent of L.

In any practical application of the approximation algorithm outlined in the foregoing section to a CNF-formula G one might proceed as follows: collect all clauses with variable $a_1$ and $\sim a_1$. Call the remaining factor $G_R$. Then one has (in the notation of (28) with $l=a_1$)

(30)    $D^{(1)}(G) = G \vee R_1(G) = G(a_1,...) \vee G(\sim a_1,...)$

$\qquad = (S_1 S_2 ... S_t \vee R_1 R_2 ... R_s) G_R$

where neither the R and S nor $G_R$ depend on $a_1$. The collection procedure is polynomial and the resulting formula is <u>not</u> longer than the original one in terms of symbols. <u>But</u> it is not a CNF formula anymore. If one wants to repeat the process and apply the same rules, one has to split (30) into two CNF formulas and apply the procedure to each. Now in effect the formula length has doubled (nearly) and one encounters the exponential behaviour typical of NP-problems. Simplifications might arise from the S and R factors, however. All of them are shorter than the clauses one started with because they do not contain $a_1$ anymore. If an S or R is reduced to a single varible l, the application of (27b) can eliminate several clauses in one stroke. From this consideration it becomes clear that an effective algorithm will involve a clever choice of consecutive variables.

*Conclusion.*

Two new formal tools to deal with propositional calculus and the problem of satisfiability were discussed; namely the notion of logical primes [5] and the metaformula. It was shown that each equivalence class of Boolean formulas has a unique representation as a product of logical primes. Therefore the satisfiability of a formula can be formulated as a problem of prime factorization.

The notion of the metavariable or metaformula enables one to formulate well known procedures for determining satisfiability in a systematic manner. A simple program was formulated which checks for SAT in n (number of basic variables) linear steps. Nonetheless the procedure cannot do the job in polynomial time because the length of the formula to be checked in each step basically doubles. Steps to optimize the procedure by proper length reductions were indicated.

*Appendix*

The proofs for theorems (TI) to (TV) are straightforward and only basic ideas will be sketched here.

Proof of TI: For N=1 V consists only of the trivial elements $0$ and $1$. Thus we assume $|V|>2$. For some nontrivial s define $K_s=\{a|axs=0\}$. Obviously $\sim s$ and $0 \in K_s$. Analogously for $K_{\sim s}$. It is easy to show that $K_s$ and $K_{\sim s}$ are subgroups of V with respect to $+$, and both have only $0$ in common. Thus each $a \in V$ has a unique decomposition $a=u+v$ where $u \in K_s$ and $v \in K_{\sim s}$. Let $|K_s|=N_s$, and $|K_{\sim s}|=N_{\sim s}$. Next we count elements which do not belong to $K_s$ or $K_{\sim s}$. Define:

$E_{Ks}(u_0) = \{u_0+v| v \in K_{\sim s}\setminus 0\}$ with $u_0 \in K_s$. $|E_{Ks}(u_0)| = N_{\sim s}-1$ from the definition. Next one shows that $E_{Ks}(a)$ and $E_{Ks}(b)$ have no elements in common unless $a=b$. Thus $|V|= N_s-1+ N_{\sim s} +|\Sigma_u E_{Ks}(u)| = N_s-1+ N_{\sim s} +(N_s-1)|E_{Ks}(u)|= (N_s-1)(1+N_{\sim s}-1)+ N_{\sim s}$

$= N_s N_{\sim s}$.

Since both $K_s$ and $K_{\sim s}$ are subfields of V (with neutral elements $\sim s$ and s with respect to x) one can apply the same line of argument to each of them until one reaches the trivial field $V_0 = \{0, \tilde{1}\}$ which has $|V_0|=2$. Thus both $N_s$ and $N_{\sim s}$, and therefore $|V|$ is a power of 2.

Next the proof of (TII) can proceed via induction over $N = ld(|V|)$.

Again one considers the subfields $K_s$ and $K_{\sim s}$ of a V with $|V| = 2^{N+1}$ and their sets of primes $p_j$ and $q_j$ which exist by assumption. Then one shows that all $p_j + s$ are primes in V, and $q_j + \sim s$ dto. Furthermore one can show that no two of these primes of V or their negations coincide, and, secondly, that any possible prime of V is necessarily one of them. Thus the $p_j + s$ and $q_j + \sim s$ constitute the set of primes of V, and their number is by assumption $ld(N_s) + ld(N_{\sim s}) = N+1$.

The fact that different negated $p_k$ are orthogonal, equ. (8), is proven as follows:

For $i \neq j$ $p_j x \sim p_i \, \varepsilon \, K_{pi}$ by definition of K. But since $p_i$ is prime, $K_{pi} = \{0, \sim p_i\}$. Thus either $p_j x \sim p_i = 0$ which implies (because also $p_j$ is prime) that $\sim p_i$ is either $0$ or equal to $\sim p_j$ both in contradiction to assumptions, therefore : or $p_j x \sim p_i = \sim p_i$. Which is equivalent to the claim.

Along the same line of thought - considering $K_s$ and $K_{\sim s}$ for s=some prime element of V - it can be proven that each element of V has a unique decomposition into primes, equ. (7) or (6).

Proof of (TIV).

First note that both functions f(x,y) and g(x,y) in equ. (9) can take values 0 or 1 only, and they are symmetric because of the commutativity of the operations x and +. Then from (1) and (9) setting T(a)=0 or 1 respectively one gets

$0 = g(0,T(\mathbf{1})) = g(T(\mathbf{1}),0)$ and $1 = g(1,T(\mathbf{1})) = g(T(\mathbf{1}),1)$ and

$T(\mathbf{0}) = g(1, T(\mathbf{0})) = g(0, T(\mathbf{0}))$ from ax $\mathbf{0}= \mathbf{0}$.

If one chooses $T(\mathbf{0}) = 0$ then $T(\mathbf{1}) = 0$ leads to a contradiction, as well as setting both values equal to 1. One is left with the choice

(A)  $T(\mathbf{0}) = 0$ and $T(\mathbf{1}) = 1$

(B)  $T(\mathbf{0}) = 1$ and $T(\mathbf{1}) = 0$

We adopt choice (A) in the following. As a consequence

$0=g(0,1) = g(1,0) = g(0,0)$ and $1 = g(1,1)$ and, from (1) for +

$0=f(0,0)=f(1,1)$ and $1=f(1,0)=f(0,1)$.

Let T be fixed. Because of (8): $0=g(T(\tilde{p}),T(\tilde{q}))$ for different p,q. Thus either $T(\tilde{p})=T(\tilde{q})=0$ or the two assignments have different value. If $T(\tilde{p}_k)=0$ for all k, one gets a contradiction to $\mathbf{1}=\Sigma_k \tilde{p}_k$ and $0=f(0,0)$. Thus at least for one k $T(\tilde{p}_k)=1$. But then for all other j $T(\tilde{p}_j)=0$ because of $0=g(0,1)$ and the orthogonality relation (8). Thus for each T there is exactly one $\tilde{p}_k$ with truth assignment 1, and all other $\tilde{p}$ giving 0. Now consider two different maps T, T' with $T(\tilde{p}_k)=1$ and $T'(\tilde{p}_l)=1$. Then k and l must be different, otherwise the two maps would coincide. Repeating this argument with a third T" and so on leads to the conclusion that there are exactly as many allowed maps as there are primes. We can label the maps as we would like to, so the most natural choice is equ. (10).

As for theorem V, the easiest way to prove the existence of n=ld(N) $a_k$ is to construct them from the uniquely defined primes:

$a_r = \Sigma_i \Sigma_s \Sigma_l \sim p_i \, \delta(i, s+2^k l)$

where $\delta$ is the Kronecker $\delta$ and the s and l sums run from $2^{k-1}$ to $2^k-1$ and from 0 to $2^{n-k}-1$ respectively. Constructing them inductively is more instructive because one encounters choices which lead to different sets of $a_k$. The seemingly complicated formula above is obsolete once one uses the binary representation of all quantities which is given by the bijection F $\longleftrightarrow$ $T_{N-1}(F) ... T_i(F) ... T_0(F)$ for any F. In particular the $a_i$ take the simple form:

$a_1 = $ ....1010101010101010

$a_2 = $ ....1100110011001100

$a_3 = $ ....1111000011110000

and so on.

References.